\documentclass{gen-p-l}
\usepackage{amsmath}
\usepackage{amsthm}
\usepackage{amscd}
\usepackage{amssymb}

\input{epsf.tex}

\numberwithin{equation}{section}

\theoremstyle{definition}

\newtheorem*{Remark}{Remark}

\theoremstyle{plain}
\newtheorem{theorem}{Theorem}

\newtheorem{corollary}[theorem]{Corollary}

\newtheorem{conjecture}[theorem]{Conjecture}
\newtheorem{question}[theorem]{Question}
\newtheorem{prob}[theorem]{Problem}

\newcounter{remarks}
\newenvironment{remarks}%
{
     \begin{list}{\arabic{remarks}. }{\usecounter{remarks}%
          \setlength{\leftmargin}{0in}%
          \setlength{\rightmargin}{0in}%
          \setlength{\labelsep}{0pt}%
          \setlength{\labelwidth}{0pt}%
          \setlength{\listparindent}{0pt}%
     }
}
{
\end{list}
}

\DeclareMathOperator{\Aff}{Aff}

\DeclareMathOperator{\vcd}{vcd}

\DeclareMathOperator{\SL}{SL}
\DeclareMathOperator{\GL}{GL}

\DeclareMathOperator{\QC}{QC}
\DeclareMathOperator{\BS}{BS}
\DeclareMathOperator{\Isom}{Isom}

\DeclareMathOperator\Bilip{Bilip}

\DeclareMathOperator\SO{SO}
\DeclareMathOperator\SU{SU}
\DeclareMathOperator\wreath{\wr}
\DeclareMathOperator\deggrowth{deg}
\DeclareMathOperator\graded{gr}

\newcommand\R{{\mathbf R}}

\newcommand\Q{{\mathbf Q}}

\newcommand\C{{\mathbf C}}

\newcommand\Z{{\mathbf Z}}

\newcommand{\from}{\colon}

\newcommand\suchthat{\bigm|}

\newcommand\absvalue[1]{\left| #1 \right|}
\newcommand\abs[1]{\absvalue{#1}}

\newcommand\semidirect{\rtimes}

\DeclareMathOperator\QI{QI}
\newcommand\cross{\times}

\newcommand\Haus{{\mathcal H}}

\renewcommand\H{{\mathbf H}}

\newcommand\N{{\mathbf N}}

\newcommand\weakto\looparrowright
\newcommand\weakfrom\looparrowleft

\newcommand\ds\displaystyle

\newcommand\solv{\ensuremath{\text{\scshape solv}}}
\newcommand\Solv\solv

\begin{document}

\title{Problems on the geometry of \\ finitely generated solvable groups}

\author{Benson Farb}

\curraddr{Department of Mathematics\\ University of Chicago\\ 
5734 University Ave.\\ Chicago, Il 60637}

\email{farb@math.uchicago.edu}

\thanks{Supported in part by NSF grant DMS 9704640 and by a Sloan
Foundation fellowship.}

\author{Lee Mosher}

\curraddr{Department of Mathematics\\ Rutgers University\\
Newark, NJ 07102}

\email{mosher@andromeda.rutgers.edu}

\thanks{Supported in part by NSF grant DMS 9504946.}

\thanks{Appears in ``Crystallographic Groups and their 
   Generalizations (Kortrijk, 1999)'', (ed. P. Igodt et. al.),
121--134, Cont. Math. 262, Amer. Math. Soc. (2000).}

\begin{abstract}
A survey of problems, conjectures, and theorems about quasi-isometric
classification and rigidity for finitely generated solvable groups.
\end{abstract}


\maketitle
\tableofcontents

\section{Introduction}

Our story begins with a theorem of Gromov, proved in 1980.

\begin{theorem}[Gromov's Polynomial Growth Theorem \cite{Gr1}]
\label{theorem:polynomial:growth}
Let $G$ be any finitely generated group.  If $G$ has polynomial growth then  
$G$ is virtually nilpotent, i.e.\ $G$ has a finite index nilpotent subgroup. 
\end{theorem}

Gromov's theorem inspired the more general problem (see, e.g. 
\cite{GH1,BW,Gr1,Gr2}) 
of understanding to what extent the asymptotic geometry of a
finitely generated solvable group determines its algebraic
structure. One way in which to pose this question precisely is via the
notion of quasi-isometry.

A (coarse) {\em quasi-isometry} between metric spaces 
is a map $f:X\to Y$ such that, for 
some constants $K,C,C'>0$:
\begin{enumerate}
\item 
$\frac{1}{K} d_X(x_1,x_2) - C \le d_Y(f(x_1),f(x_2)) \le K
d_X(x_1,x_2) + C$ for all $x_1,x_2\in X$.

\item The $C'$-neighborhood of $f(X)$ is all of $Y$.
\end{enumerate}
\medskip
$X$ and $Y$ are \emph{quasi-isometric} if there exists a quasi-isometry
$X \to Y$. Note that the quasi-isometry type of a metric space $X$ is
unchanged upon removal of any bounded subset of $X$; hence the term
``asymptotic''.

Quasi-isometries are the natural maps to study when one is interested in 
the geometry of a group.  In particular:

\medskip

\begin{remarks}
\item 
The word metric on any f.g.\ group is unique up to quasi-isometry.

\item 
Any injective homomorphism with finite index image is a quasi-isometry,
as is any surjective homomorphism with finite kernel. The
equivalence relation generated by these two types of  maps can be
described more compactly: two groups
$G,H$ are equivalent in this manner if and only if they are \emph{weakly
commensurable}, which means that there exists a group $Q$ and
homomorphisms $Q\to G$, $Q \to H$ each having finite kernel and finite
index image (proof: show that ``weakly commensurable'' is in fact an
equivalence relation). This weakens the usual notion of {\em
commensurability}, i.e.\ when $G$ and $H$ have isomorphic finite index
subgroups. Weakly commensurable groups are clearly
quasi-isometric.  

\item 
Any two cocompact, discrete subgroups of a Lie group are
quasi-isometric. There are cocompact discrete subgroups
of the same Lie group which are not weakly commensurable, for example
arithmetic lattices are not weakly commensurable to non-arithmetic ones. 

\end{remarks}

\medskip

The Polynomial Growth Theorem was an important motivation for Gromov
when he initiated in \cite{Gr2,Gr3} the problem of classifying
finitely-generated groups up to quasi-isometry.

Theorem \ref{theorem:polynomial:growth}, together with the fact that
nilpotent groups have polynomial growth (see \S\ref{section:nilpotent}
below), implies that the property of
being nilpotent is actually an asymptotic property of groups.  More
precisely, the class of nilpotent groups is {\em quasi-isometrically
rigid}: any finitely-generated group quasi-isometric to a nilpotent
group is weakly commensurable to some nilpotent group
\footnote{In fact such a group must have a finite-index nilpotent
subgroup.  In general, weak commensurability is the most one can hope
for in quasi-isometric rigidity problems.}.  Sometimes this is
expressed by saying that the property of being nilpotent is a {\em
geometric property}, i.e.\ it is a quasi-isometry invariant (up to 
weak commensurability).  The natural question then becomes:

\begin{question}[Rigidity question]
\label{question:rigidity}
Which subclasses of f.g.\ solvable groups 
are quasi-isometrically rigid? For example, are polycyclic groups 
quasi-isometrically rigid? metabelian groups? nilpotent-by-cyclic groups?
\end{question}

In other words, which of these algebraic properties of a group are
actually geometric, and are determined by apparently cruder asymptotic
information?  A.\ Dioubina \cite{Di} has recently found examples which
show that the class of finitely generated solvable groups is {\em not}
quasi-isometrically rigid (see \S\ref{section:dioubina}).
On the other hand, at least some subclasses of solvable groups are
indeed rigid (see \S\ref{section:abc}).

Along with Question \ref{question:rigidity} comes the finer
classification problem:

\begin{prob}[Classification problem]
\label{problem:classification}
Classify f.g.\ solvable (resp. nilpotent, polycyclic, 
metabelian, nilpotent-by-cyclic, etc.) groups up to quasi-isometry.
\end{prob}

As we shall see, the classification problem is usually much more
delicate than the rigidity problem; indeed the quasi-isometry
classification of finitely-generated nilpotent groups remains one of the
major open problems in the field.  We discuss this in greater detail in
\S\ref{section:nilpotent}.

The corresponding rigidity and classification problems for irreducible
lattices in semisimple Lie groups have been completely solved.  This is
a vast result due to many people and, among other things, it generalizes
and strengthens the Mostow Rigidity Theorem.  We refer the reader to
\cite{Fa} for a survey of this work.

In contrast, results for finitely generated solvable groups have remained 
more elusive.  There are several reasons for this:
\begin{remarks}

\item Finitely generated solvable groups are  defined algebraically, and
so they do not always come equipped with an obvious or well-studied
geometric model (see, e.g., item 3 below).

\item Dioubina's examples show not only that the class 
of finitely-generated solvable groups is not quasi-isometrically rigid;
they also show (see \S\ref{section:abc} below) that the answer to
Question \ref{question:rigidity} for certain subclasses of solvable 
groups (e.g.\ abelian-by-cyclic) differs in the finitely presented and
finitely generated cases.

\item There exists a finitely presented solvable group $\Gamma$ 
of derived length 3 with the property that $\Gamma$ has unsolvable word
problem (see \cite{Kh}).  Solving the word problem for a group is
equivalent to giving an algorithm to build the Cayley graph of that
group.  In this sense there are finitely presented solvable groups whose
geometry cannot be understood, at least by a Turing machine.

\item Solvable groups are much less rigid than irreducible lattices in
semisimple Lie groups.  This phenomenon is exhibited concretely by the
fact that many finitely generated solvable groups have
infinite-dimensional groups of self quasi-isometries\footnote{The set of self quasi-isometries of a space $X$, with
the operation of composition, becomes a group $\QI(X)$ once one mods out
by the relation $f\sim g$ if $d(f,g)<\infty$ in the $\sup$ norm.} (see
below).
\end{remarks}

\begin{prob}[Flexibility of solvable groups]
\label{problem:flexible}
For which infinite, finitely generated solvable groups $\Gamma$ is
$\QI(\Gamma)$ infinite dimensional?
\end{prob}

\noindent
In contrast, all irreducible lattices in semisimple Lie groups $G\neq
\SO(n,1), \SU(n,1)$ have countable or finite-dimensional quasi-isometry
groups.

\bigskip

At this point in time, our understanding of the geometry of
finitely-generated solvable groups is quite limited.  In
\S\ref{section:nilpotent} we discuss what is known about the
quasi-isometry classification of nilpotent groups (the rigidity being
given by Gromov's Polynomial Growth Theorem).  Beyond nilpotent groups,
the only detailed knowledge we have is for the finitely-presented,
nonpolycyclic abelian-by-cyclic groups.  We discuss this in depth in
\S\ref{section:abc}, and give a conjectural picture of the polycyclic
case in \S\ref{section:abc2}.  One of the interesting discoveries
described in these sections is a connection between finitely presented
solvable groups and the theory of dynamical systems.  This connection is
pursued very briefly in a more general context in
\S\ref{section:final}, together with some questions about issues beyond
the limits of current knowledge.

This article is meant only as a brief survey of problems, conjectures,
and theorems. It therefore contains neither an exhaustive history nor
detailed proofs; for these the reader may consult the references. It is a
pleasure to thank David Fisher, Pierre de la Harpe, Ashley Reiter,
Jennifer Taback, and the referee for their comments and corrections.

\section{Dioubina's examples}
\label{section:dioubina}

Recall that the {\em wreath product} of groups $A$ and $B$, denoted
$A\wreath B$, is the semidirect product
$(\oplus_A B) \semidirect A$, where $\oplus_A B$ is the direct sum of
copies of B indexed by elements of A, and A acts via the ``shift'',
i.e.\ the left action of $A$ on the index set $A$ via left
multiplication.  Note that if $A$ and $B$ are finitely-generated then so
is $A\wreath B$.

The main result of Dioubina \cite{Di} is that, if there is a bijective 
quasi-isometry between finitely-generated 
groups $A$ and $B$, then for any finitely-generated group $C$ the groups 
$C\wreath A$ and $C\wreath B$ are quasi-isometric.  Dioubina then
applies this theorem to the groups $A=C=\Z, B=\Z\oplus D$ where $D$ is a
finite nonsolvable group.  It is easy to construct a one-to-one
quasi-isometry between $A$ and $B$.  Hence 
$G=C\wreath A$ and $H=C\wreath B$ are quasi-isometric.  

Now $G$ is torsion-free solvable, in fact $G=\Z\wreath \Z$ is an
abelian-by-cyclic group of the form $\Z[\Z]$-by-$\Z$. On the other hand
$H$ contains $\oplus_{\Z}D$, and so is not virtually solvable, nor even
weakly commensurable with a solvable group.  Hence the class of
finitely-generated solvable groups is not quasi-isometrically rigid.

Dioubina's examples never have any finite 
presentation.  In fact if $A\wreath B$ is finitely presented then 
either $A$ or $B$ is finite (see \cite{Bau}).  
This leads to the following question.

\begin{question}
Is the class of finitely presented solvable groups quasi-isometrically
rigid?
\end{question}

Note that the property of being finitely presented is a quasi-isometry
invariant (see \cite{GH1}).  

\section{Nilpotent groups and Pansu's Theorem}
\label{section:nilpotent}

While the Polynomial Growth Theorem shows that the class of finitely
generated nilpotent groups is quasi-isometrically rigid, the following
remains an important open problem.

\begin{prob}
Classify finitely generated nilpotent groups up to quasi-isometry.
\end{prob}

The basic quasi-isometry invariants for a finitely-generated nilpotent
group $G$ are most easily computed in terms of the set $\{c_i(G)\}$ 
of ranks (over $\Q$) of the quotients $G_i/G_{i+1}$ of the lower central 
series for $G$, where $G_i$ is defined inductively by $G_1=G$ and 
$G_{i+1}=[G,G_i]$.  

One of the first quasi-isometry invariants to be studied was the growth
of a group, studied by Dixmier, Guivar'ch, Milnor, Wolf, and others (see
\cite{H}, Chapters VI-VII for a nice discussion of this, and a careful
account of the history).  The {\em growth} of $G$ is the function of
$r$ that counts the number of elements in a ball of radius $r$ in
$G$. There is an important dichotomy for solvable groups:

\begin{theorem}[\cite{Wo,Mi}]
Let $G$ be a finitely generated solvable group.  Then either $G$ has
polynomial growth and is virtually nilpotent, or $G$ has exponential 
growth and is not virtually nilpotent.
\end{theorem}

When $G$ has polynomial growth, the degree $\deggrowth(G)$ of this
polynomial is easily seen to be a quasi-isometry invariant. It is given by
the following formula, discovered around the same time by Guivar'ch
\cite{Gu} and by Bass \cite{Ba}:
$$\deggrowth(G)=\sum_{i=1}^n i \, c_i(G)
$$ 
where $n$ is the degree of nilpotency of $G$.  Another basic invariant is
that of virtual cohomological dimension $\vcd(G)$.  For groups $G$ with
finite classifying space (which is not difficult to check for torsion-free
nilpotent groups), this number was shown by Gersten \cite{Ge} and
Block-Weinberger \cite{BW} to be a quasi-isometry invariant.  On the other
hand it is easy to check that 
$$\vcd(G)=\sum_{i=1}^nc_i(G)
$$
where $n$ is the degree of nilpotency, also known as the Hirsch length, of
$G$.  As Bridson and Gersten have shown (see \cite{BG,Ge}), the above two
formulas imply that any finitely generated group $\Gamma$ which is
quasi-isometric to $\Z^n$ must have a finite index $\Z^n$ subgroup: by the
Polynomial Growth Theorem such a $\Gamma$ has a finite index nilpotent
subgroup $N$; but 
$$\deggrowth(N)=\deggrowth(\Z^n)=n=\vcd(\Z^n)=\vcd(N)$$
and so 
$$\sum_i i \, c_i(N)=\sum_i c_i(N)
$$
which can only happen if $c_i(N)=0$ for $i>1$, in which case $N$ is
abelian.  

\begin{prob}[\cite{GH2}]
Give an elementary proof (i.e. without using Gromov's Polynomial Growth
Theorem) that any finitely generated group quasi-isometric to $\Z^n$ has 
a finite index $\Z^n$ subgroup.
\end{prob}

As an exercise, the reader is invited to find nilpotent
groups $N_1,N_2$ which are not quasi-isometric but which have the same
degree of growth and the same $\vcd$.

There are many other quasi-isometry invariants for finitely-generated 
nilpotent groups $\Gamma$.  All known invariants are special cases of
the following theorem of Pansu \cite{Pa1}.  To every nilpotent group
$\Gamma$ one can associate a nilpotent Lie group $\Gamma \otimes \R$,
called the {\em Malcev completion} of $\Gamma$ (see \cite{Ma}), as well
as the associated graded Lie group $\graded(\Gamma\otimes \R)$. 

\begin{theorem}[Pansu's Theorem \cite{Pa1}]
\label{theorem:pansu}
Let $\Gamma_1,\Gamma_2$ be two finitely-generated nilpotent groups.  If 
$\Gamma_1$ is quasi-isometric to $\Gamma_2$ then
$\graded(\Gamma_1\otimes \R)$ is isomorphic to $\graded(\Gamma_2\otimes
\R)$. 
\end{theorem}

We remark that there are nilpotent groups with non-isomorphic Malcev
completions where the associated gradeds are isomorphic; the examples
are 7-dimensional and somewhat involved (see \cite{Go}, p.24, Example
2).  It is not known whether or not the Malcev completion is a
quasi-isometry invariant.

Theorem \ref{theorem:pansu} immediately implies:

\begin{corollary}
The numbers $c_i(\Gamma)$ are quasi-isometry invariants.  
\end{corollary}

In particular we recover (as special cases) that growth and
cohomological dimension are quasi-isometry invariants of $\Gamma$.

To understand Pansu's proof one must consider {\em Carnot groups}.
These are graded nilpotent Lie groups $N$ whose Lie algebra $\N$
is generated (via bracket) by elements of degree one.  Chow's Theorem 
\cite{Ch} 
states that such Lie groups $N$ have the property that the
left-invariant distribution obtained from the degree one subspace $\N^1$ of $\N$ is a {\em totally nonintegrable} distribution: any
two points $x,y\in N$ can be connected by a piecewise smooth path
$\gamma$ in $N$ for which the vector $d\gamma/dt(t)$ lies in the
distribution.  Infimizing the length of such paths between two given
points gives a metric on $N$, called the {\em Carnot Careth\'{e}odory}
metric $d_{car}$.  This metric is non-Riemannian if $N\neq \R^n$.  For
example, when $N$ is the $3$-dimensional Heisenberg group then the
metric space $(N,d_{car})$ has Hausdorff dimension $4$.

One important property of Carnot groups is that they come equipped with
a $1$-parameter family of dilations $\{\delta_t\}$, which gives 
a notion of {\em (Carnot) differentiability} (see \cite{Pa1}).  Further, the
differential $Df(x)$ of a map $f:N_1\rightarrow N_2$ between Carnot
groups $N_1,N_2$ which is (Carnot) differentiable at the point $x\in
N_1$ is actually a {\em Lie group homomorphism} $N_1\rightarrow N_2$.

\bigskip
\noindent
{\bfseries Sketch of Pansu's proof of Theorem \ref{theorem:pansu}.\/}  
If $(\Gamma,d)$ is a nilpotent group endowed with a word metric
$d$, the sequence of scaled metric spaces $\{(\Gamma,
\frac{1}{n}d)\}_{n\in
\N}$ has a limit in the sense of Gromov-Hausdorff convergence:
$$(\Gamma_{\infty},d_{\infty})=\lim_{n\rightarrow
\infty}(\Gamma,\frac{1}{n}d)
$$
(See \cite{Pa2}, and \cite{BrS} for an introduction to Gromov-Hausdorff
convergence). It was already known, using ultralimits, that some
subsequence converges \cite{Gr1,DW}. Pansu's proof not only gives
convergence on the nose, but it yields some additional important features
of the limit metric space $(\Gamma_{\infty},d_{\infty})$: 
\begin{itemize}
\item (Identifying limit) 
It is isometric to the Carnot group $\graded(\Gamma\otimes\R)$
endowed with the Carnot metric $d_{car}$.

\item (Functoriality) Any quasi-isometry $f:\Gamma_1\rightarrow
\Gamma_2$ bewteen finitely-generated nilpotent groups induces a
{\em bilipschitz homeomorphism} 
$$\hat{f}:(\graded(\Gamma_1\otimes\R),d_{car})\rightarrow
(\graded(\Gamma_2\otimes\R),d_{car})$$  
\end{itemize}

Note that functoriality follows immediately once we know the limit
exists: the point is that if $f:(\Gamma_1,d_1)\rightarrow
(\Gamma_2,d_2)$ is a $(K,C)$ quasi-isometry of word metrics, then for 
each $n$ the map 
$$f:(\Gamma_1,\frac{1}{n}d_1)\rightarrow (\Gamma_2,\frac{1}{n}d_2)$$ is a
$(K,C/n)$ quasi-isometry, hence the induced map
$$\hat{f}:((\Gamma_1)_\infty,d_{car})\rightarrow
((\Gamma_2)_\infty,d_{car})$$ is a $(K,0)$ quasi-isometry, i.e.\ is a
bilipschitz homeomorphism.

Given a quasi-isometry $f:\Gamma_1\rightarrow \Gamma_2$, we thus
have an induced bilipschitz homeomorphism
$\hat{f}:\graded(\Gamma_1\otimes \R)\rightarrow \graded(\Gamma_2\otimes
\R)$ between Carnot groups endowed with 
Carnot-Careth\'{e}odory metrics.  Pansu then
proves a regularity theorem, generalizing the Rademacher-Stepanov
Theorem for $\R^n$.  This general regularity theorem 
states that a bilipschitz homeomorphism of Carnot groups (endowed with 
Carnot-Careth\'{e}odory metrics) is differentiable almost everywhere.  
Since the differential
$D_x\hat{f}$ is actually a group homomorphism, we know that for almost
every point $x$ the differential 
$D_x\hat{f}:\graded(\Gamma_1\otimes \R)\rightarrow \graded(\Gamma_2\otimes
\R)$ is an isomorphism.

\section{Abelian-by-cyclic groups: nonpolycyclic case}
\label{section:abc}

The first progress on Question \ref{question:rigidity} and Problem
\ref{problem:classification} in the non-(virtually)-nilpotent case was made in
\cite{FM1} and \cite{FM2}.  These papers proved classification and
rigidity for the simplest class of non-nilpotent solvable groups: the
{\em solvable Baumslag-Solitar groups} $$\BS(1,n)=<a,b: aba^{-1}=b^n>$$

These groups are part of the much broader class of abelian-by-cyclic 
groups.  A group $\Gamma$ is {\em abelian-by-cyclic}
if there is an exact sequence $$1\rightarrow A\rightarrow
\Gamma\rightarrow Z\rightarrow 1 $$ where $A$ is an abelian group and
$Z$ is an infinite cyclic group. If $\Gamma$ is finitely generated, then
$A$ is a finitely generated module over the group ring $\Z[Z]$, although
$A$ need not be finitely generated as a group.

By a result of Bieri and Strebel \cite{BS1}, the
class of finitely presented, torsion-free, 
abelian-by-cyclic groups may be described
in another way. Consider an $n \cross n$ matrix $M$ with integral
entries and $\det M \ne 0$. Let $\Gamma_M$ be the ascending HNN
extension of $\Z^n$ given by the monomorphism $\phi_M$ with matrix
$M$. Then $\Gamma_M$ has a finite presentation
$$
\Gamma_M=\langle t,a_1,\ldots ,a_n 
\suchthat [a_i,a_j]=1, ta_it^{-1}=\phi_M(a_i),
i,j=1,\ldots,n\rangle $$ where $\phi_M(a_i)$ is the word
$a_1^{m_1}\cdots a_n^{m_n}$ and the vector $(m_1,\ldots ,m_n)$ is the
$i^{\text{th}}$ column of the matrix $M$. Such groups $\Gamma_M$ are
precisely the class of finitely presented, torsion-free,
abelian-by-cyclic groups (see \cite{BS1} for a proof involving a
precursor of the Bieri-Neumann-Strebel invariant, or
\cite{FM2} for a proof using trees). The group $\Gamma_M$ is
polycyclic if and only if $\abs{\det M}=1$ (see \cite{BS2}). 

The results of \cite{FM1} and \cite{FM2} are generalized in 
\cite{FM3}, which gives the complete classification of 
the finitely presented, nonpolycyclic abelian-by-cyclic groups among 
all f.g.\ groups, as given by the following two theorems.  

The first theorem in \cite{FM3} gives a classification of all
finitely-presented, nonpolycyclic, abelian-by-cyclic groups up to
quasi-isometry. It is easy to see that any such group has a
torsion-free subgroup of finite index, so is commensurable (hence
quasi-isometric) to some $\Gamma_M$. The classification of these groups
is actually quite delicate---the standard quasi-isometry
invariants (ends, growth, isoperimetric inequalities, etc.) do not
distinguish any of these groups from each other, except that the size of
the matrix $M$ can be detected by large scale cohomological invariants of
$\Gamma_M$.

Given $M\in \GL(n,\R)$, the \emph{absolute Jordan form} of $M$ is the
matrix obtained from the Jordan form for $M$ over $\C$ by replacing each
diagonal entry with its absolute value, and rearranging the Jordan blocks
in some canonical order. 

\begin{theorem}[Nonpolycyclic, abelian-by-cyclic groups: Classification]
\label{theorem:classification}
Let $M_1$ and $M_2$ be integral matrices with $\abs{\det M_i}>1$ for
$i=1,2$. Then $\Gamma_{M_1}$ is quasi-isometric to 
$\Gamma_{M_2}$ if and only if there are positive integers $r_1,r_2$ such 
that $M_1^{r_1}$ and $M_2^{r_2}$ have the same absolute Jordan form.
\end{theorem}

\begin{Remark} 
Theorem \ref{theorem:classification} generalizes the main result of
\cite{FM1}, which is the case when $M_1, M_2$ are positive
$1\times 1$ matrices; in that case the result of \cite{FM1}
says even more, namely that $\Gamma_{M_1}$ and $\Gamma_{M_2}$ are
quasi-isometric if and only if they are commensurable. When $n \ge 2$,
however, it's not hard to find $n\cross n$ matrices $M_1, M_2$ such that
$\Gamma_{M_1}, \Gamma_{M_2}$ are quasi-isometric but not commensurable.
Polycyclic examples are given in \cite{BG}; similar ideas 
can  be used to produce nonpolycyclic examples.
\end{Remark}

The following theorem shows that the algebraic property of being
a finitely presented, nonpolycyclic, abelian-by-cyclic group is in fact
a geometric property.

\begin{theorem}[Nonpolycyclic, abelian-by-cyclic groups: Rigidity]
\label{theorem:rigidity}
Let $\Gamma=\Gamma_M$ be a finitely presented abelian-by-cyclic group,
determined by an $n \cross n$ integer matrix $M$ with $\abs{\det M}>1$.
Let $G$ be any finitely generated group quasi-isometric to $\Gamma$. Then
there is a finite normal subgroup $N \subset G$ such that $G/N$ is 
commensurable to $\Gamma_N$, for some $n \cross n$ integer
matrix $N$ with $\abs{\det N}>1$.
\end{theorem}

\begin{Remark}
Theorem \ref{theorem:rigidity} generalizes the main result of 
\cite{FM2}, which covers the case when $M$ is a positive
$1\times 1$ matrix. The $1\times 1$ case is given a new proof in
\cite{MSW}, which is adapted in \cite{FM3} to prove Theorem
\ref{theorem:rigidity}.
\end{Remark}

\begin{Remark}
The ``finitely presented'' hypothesis in Theorem \ref{theorem:rigidity}
cannot be weakened to ``finitely generated'', since Diuobina's example
(discussed in \S\ref{section:dioubina}) is abelian-by-cyclic, namely 
$\Z[\Z]$-by-$\Z$. 
\end{Remark}

One new discovery in \cite{FM3} is that there is a strong connection
between the geometry of solvable groups and the theory of dynamical
systems.  Assuming here for simplicity that the matrix $M$ lies on a
$1$-parameter subgroup $M^t$ in $\GL(n,\R)$, let $G_M$ be the
semi-direct product $\R^n\semidirect_M \R$, where $\R$ acts on $\R^n$ by
the $1$-parameter subgroup $M^t$. We endow the solvable Lie group $G_M$
with a left-invariant metric.  The group $G_M$ admits a {\em vertical
flow}: 
$$\Psi_s(x,t)=(x,t+s)$$ 
There is a natural {\em horizontal
foliation} of $G_M$ whose leaves are the level sets $P_t=\R^n\times
\{t\}$ of time.  A quasi-isometry $f:G_M\to G_N$ is {\em horizontal
respecting} if it coarsely permutes the leaves of this foliation; that
is, if there is a constant $C\geq 0$ so that
$$d_\Haus(f(P_t),P_{h(t)})\leq C$$ where $d_\Haus$ denotes Hausdorff
distance and $h:\R\rightarrow
\R$ is some function, which we think of as a {\em time change} between
the flows.  

A key technical result of \cite{FM3} is the phenomenon of {\em time
rigidity}: the time change $h$ must actually be {\em affine}, so taking a
real power of $M$ allows one to assume $h(t)=t$.  

It is then shown that ``quasi-isometries remember the dynamics''. That
is, $f$ coarsely respects several foliations arising from the partially
hyperbolic dynamics of the flow $\Psi$, starting with the weak stable,
weak unstable, and center-leaf foliations. By keeping track of different
exponential and polynomial divergence properties of the action of $\Psi$
on tangent vectors, the weak stable and weak unstable foliations are
decomposed into flags of foliations. Using time rigidity and an
inductive argument it is shown that these flags are coarsely respected by
$f$ as well. Relating the flags of foliations to the Jordan Decomposition
then completes the proof of:

\begin{theorem}[Horizontal respecting quasi-isometries]
\label{theorem:horizontal}
If there is a horizontal-respecting quasi-isometry $f:G_M\to G_N$ 
then there exist nonzero $a,b\in \R$ so that $M^{a}$ and $M^{b}$ 
have the same absolute Jordan form.
\end{theorem}

The ``nonpolycyclic'' hypothesis (i.e.\ $\abs{\det M}>1$) in 
Theorem \ref{theorem:classification} is used in two ways.  
First, the group $\Gamma_M$ has a model space which is 
topologically a product of $\R^n$ and a regular tree of valence 
$\abs{\det M}+1$, and when this valence is greater than $2$ we can use
coarse  algebraic topology (as developed in \cite{FS}, \cite{EF}, and
\cite{FM1})  to show that any quasi-isometry $\Gamma_M\rightarrow
\Gamma_N$  induces a quasi-isometry $G_M\rightarrow G_N$ satisfying the 
hypothesis of Theorem \ref{theorem:horizontal}.  Second, we are able to pick 
off {\em integral} $a,b$ by developing a ``boundary theory'' for 
$\Gamma_M$; in case $\abs{\det M}>1$ this boundary is 
a self-similar Cantor set whose bilipschitz geometry detects the 
primitive integral power of $\abs{\det M}$ by Cooper's Theorem
\cite{FM1}, finishing the proof of Theorem \ref{theorem:classification}.

\begin{prob}[Nilpotent-by-cyclic groups]
\label{problem:extend}
Extend Theorem \ref{theorem:classification} and Theorem 
\ref{theorem:rigidity} to the class of finitely-presented
nilpotent-by-cyclic groups. 
\end{prob}

Of course, as the classification of finitely-generated nilpotent groups
is still open, Problem \ref{problem:extend} is meant in the sense of
reducing the nilpotent-by-cyclic case to the nilpotent case, together
with another invariant.  This second invariant for a nilpotent-by-cyclic 
group $G$ will perhaps be the absolute
Jordan form of the matrix which is given by the action of the generator
of the cyclic quotient of $G$ on the nilpotent kernel of $G$.

\section{Abelian-by-cyclic groups: polycyclic case}
\label{section:abc2}

The polycyclic, abelian-by-cyclic groups 
are those $\Gamma_M$ for which $\abs{\det M}=1$, 
so that $\Gamma_M$ is cocompact 
and discrete in $G_M$, hence quasi-isometric to $G_M$.  In this case 
the proof of Theorem \ref{theorem:classification} outlined above 
breaks down, but this is so in part because the answer is quite 
different: the quasi-isometry classes of polycyclic $\Gamma_M$ are 
much coarser than in the nonpolycyclic case, as the former are 
(conjecturally) determined by 
the absolute Jordan form up to {\em real}, as opposed to integral, 
powers.  The key conjecture is:

\begin{conjecture}[Horizontal preserving]
\label{conjecture:poly1}
Suppose that $\abs{\det M}, \abs{\det N}=1$, and that $M$ and $N$ have no
eigenvalues on the unit circle.  Then every quasi-isometry of $G_M \to
G_N$ is horizontal-respecting.
\end{conjecture}

The general (with arbitrary eigenvalues) 
case of Conjecture \ref{conjecture:poly1}, which is slightly 
more complicated to state, together with Theorem \ref{theorem:horizontal} 
easily implies: 

\begin{conjecture}[Classification of polycyclic, abelian-by-cyclic groups]
Suppose that $\abs{\det M}, \abs{\det N}=1$.  Then $\Gamma_M$ 
is quasi-isometric to $\Gamma_N$ if and only if there exist 
nonzero $a,b\in \R$ 
so that $M^a$ and $N^b$ have the same absolute Jordan form. 
\end{conjecture}

Here by $M^a$ we mean $\phi(a)$, where $\phi:\R\to \GL(n,\R)$ is a
$1$-parameter subgroup with $\phi(1)=M$ (we are assuming that $M$ lies
on such a subgroup, which can be assumed after squaring $M$).  

Now let us concentrate on the simplest non-nilpotent example, which 
is also one of the 
central open problems in the field.  The 3-dimensional geometry 
\Solv\  is the Lie group 
$G_M$ where $M\in \SL(2,\Z)$ is any matrix with 2 distinct real 
eigenvalues (up to scaling, it doesn't matter which such $M$ is chosen).  

\begin{conjecture}[Rigidity of \Solv]
\label{conjecture:solv}
The $3$-dimensional Lie group \Solv\ is quasi-isometrically rigid: any f.g.\ 
group $G$ quasi-isometric to \Solv\ is weakly commensurable with a
cocompact,  discrete subgroup of \Solv. 
\end{conjecture}

There is a natural boundary for \Solv\ which decomposes into two pieces 
$\partial^s \Solv$ and $\partial^u \Solv$; these are the leaf spaces 
of the weak stable and weak unstable foliations, respectively, of 
the vertical flow on $\Solv$, and are both homeomorphic to $\R$.  

The isometry group $\Isom(\Solv)$ acts on the pair $(\partial^s \Solv, 
\partial^u \Solv)$ affinely and induces a faithful representation 
$\Solv \rightarrow \Aff(\R)\times \Aff(\R)$ whose image consists of the 
pairs 
$$(ax+b,a^{-1}x+c),\ \  a\in \R^+, b,c\in \R$$

Just as quasi-isometries of hyperbolic space $\H^n$ are characterized by 
their quasiconformal action on $\partial \H^n$ (a fact proved by Mostow), 
giving the formula $\QI(H^n)=\QC(\partial \H^n)$, we conjecture:

\begin{conjecture}[QI group of \Solv]
\label{conjecture:qigroup}
$$\QI(\Solv)=(\Bilip(\R)\times \Bilip(\R))\semidirect \ \Z/2$$
where $\Bilip(\R)$ denotes the group of bilipschitz homeomorphisms 
of $\R$, and $\Z/2$ acts by switching factors.
\end{conjecture}

There is evidence for Conjecture \ref{conjecture:qigroup}: the direction
$\supseteq$ is not hard to check (see \cite{FS}), and the analogous
theorem $\QI(\BS(1,n))=\Bilip (\R)\times \Bilip (\Q_n)$ was proved in
\cite{FM1}.  By using convergence groups techniques and a theorem of
Hinkkanen on uniformly quasisymmetric groups (see \cite{FM2}), we have
been able to show:

\begin{center}Conjecture \ref{conjecture:poly1} (in the $2\times 2$ case)
$\Rightarrow$ Conjecture \ref{conjecture:qigroup} $\Rightarrow$ Conjecture 
\ref{conjecture:solv}
\end{center}

Here is a restatement of Conjecture \ref{conjecture:poly1} in the $2
\times 2$ case:

\begin{conjecture}
Every quasi-isometry $f \from \Solv \to \Solv$ is horizontal respecting.
\label{conjecture:solvhorizontal}
\end{conjecture}

Here is one way \emph{not} to prove Conjecture
\ref{conjecture:solvhorizontal}. 

One of the major steps of \cite{FM1} in studying $\BS(1,n)$ was to
construct a model space $X_n$ for the group $\BS(1,n)$, study the
collection of isometrically embedded hyperbolic planes in $X_n$, and
prove that for any quasi-isometric embedding of the hyperbolic plane into
$X_n$, the image has finite Hausdorff distance from some isometrically
embedded hyperbolic plane. 

\begin{figure}
\centerline{\epsfbox{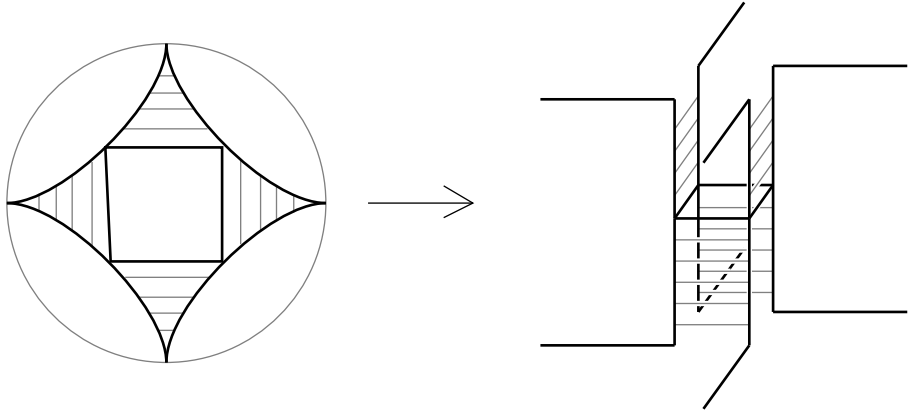}}
\caption{To map $\H^2$ quasi-isometrically into \Solv,
take a regular ideal quadrilateral in $\H^2$ divided into a regular
inscribed square and four triangles each with one ideal vertex. Map the
square isometrically to a square in the
$xy$-plane with sides parallel to the axes. Map each triangle
isometrically to a triangle in either an $xt$-plane or a $yt$-plane,
alternating around the four sides of the square. Finally, map each
complementary half-plane of the quadrilateral isometrically to a
half-plane of either an $xt$-plane or a $yt$-plane.}
\label{FigureWindvane}
\end{figure}

However, \Solv\ has quasi-isometrically embedded hyperbolic planes which
are \emph{not} Hausdorff close to isometrically embedded ones. The natural
left invariant metric on \Solv\ has the form
$$e^{2t} dx^2 + e^{-2t} dy^2 + dt^2
$$
from which it follows that the $xt$-planes and $yt$-planes are the
isometrically embedded hyperbolic planes. But none of these planes is
Hausdorff close to the set
$$\{(x,y,t)\in\Solv:  x \ge 0\ \mbox{and}\ y=0\}\cup
\{(x,y,t)\in\Solv:  y \ge 0\ \mbox{and}\ x=0\}$$
which is a quasi-isometrically embedded hyperbolic plane. An even stranger
example is shown in Figure~\ref{FigureWindvane}.

These strange quasi-isometric embeddings from $\H^2$ to \Solv\ do
share an interesting property with the standard isometric embeddings,
which may point the way to understanding quasi-isometric rigidity of
\Solv. We say that a quasi-isometric embedding
$\phi\from\H^2\to\Solv$ is \emph{$A$-quasivertical} if for each $x \in
\H^2$ there exists a vertical line $\ell \subset \Solv$ such that
$\phi(x)$ is contained in the $A$-neighborhood of $\ell$, and $\ell$ is
contained in the $A$-neighborhood of $\phi(\H^2)$.

In order to study \Solv, it therefore becomes important to understand
whether every quasi-isometrically embedded hyperbolic plane is
quasi-vertical. Specifically:

\begin{prob} Show that for all $K,C$ there exists $A$ such that each
$K,C$-quasi-isometrically embedded hyperbolic plane in \Solv\ is
$A$-quasivertical.
\label{ProblemQuasivertical}
\end{prob}

Arguing by contradiction, if Problem \ref{ProblemQuasivertical} were
impossible, fixing $K,C$ and taking a sequence of examples whose
quasi-vertical constant $A$ goes to infinity, one can pass to a
subsequence and take a renormalized limit to produce a quasi-isometric
embedding $\H^2 \to
\solv$ whose image is entirely contained in the upper half $t \ge 0$ of
\solv. But we conjecture that this is impossible:

\begin{conjecture}
There does not exist a quasi-isometric embedding $\H^2 \to \solv$ whose
image is entirely contained in the upper half space $\solv^+=\{(x,y,t)\in
\solv: t\geq 0\}$.
\end{conjecture}

\section{Next steps}
\label{section:final}

While we have already seen that there is a somewhat fine classification
of finitely presented, nonpolycyclic abelian-by-cyclic groups up to
quasi-isometry, this class of groups is but a very special class of
finitely generated solvable groups.  We have only exposed the tip of a
huge iceberg.  An important next layer is:

\begin{prob}[Metabelian groups]
Classify the finitely presented (nonpolycyclic) metabelian groups up to
quasi-isometry.  
\end{prob}

The first step in attacking this problem is to find a workable method of 
describing the geometry of the natural geometric model of such groups $G$.
Such a model should fiber over $\R^n$, where $n$ is the rank of the
maximal abelian quotient of $G$; inverse images under this projection 
of (translates of) the coordinate axes should be copies of the geometric 
models of abelian-by-cyclic groups.

\bigskip
\noindent
{\bfseries Polycyclic versus nonpolycyclic.\/}
We've seen the difference, at least in
the abelian-by-cyclic case, between polycyclic and nonpolycyclic groups.
Geometrically these two classes can be distinguished by the trees on
which they act: such trees are lines in the former case and
infinite-ended in the latter.  It is this branching
behavior which should combine with coarse topology to make the
nonpolycyclic groups more amenable to attack.  

Note that a (virtually) polycyclic group is never quasi-isometric to a
(virtually) nonpolycyclic solvable group.  This follows from the theorem
of Bieri that polycylic groups are precisely those solvable groups
satisfying Poincare duality, together with the quasi-isometric
invariance of the latter property (proved by Gersten \cite{Ge} and
Block-Weinberger
\cite{BW}). 

\bigskip
\noindent {\bfseries Solvable groups as dynamical systems.\/}
The connection of nilpotent groups with dynamical systems was made
evident in \cite{Gr1}, where Gromov's Polynomial Growth Theorem was the
final ingredient, combining with earlier work of Franks and Shub
\cite{Sh}, in the positive solution of the {\em Expanding Maps
Conjecture}: every locally distance expanding map on a closed manifold
$M$ is topologically conjugate to an expanding algebraic endomorphism on
an infranil manifold (see \cite{Gr1}).

In \S\ref{section:abc} and \S\ref{section:abc2} we saw in another way how
invariants from dynamics give quasi-isometry invariants for
abelian-by-cyclic groups.  This should be no big surprise: after all, a
finitely presented abelian-by-cyclic group is describable up to
commensurability as an ascending HNN extension $\Gamma_M$ over a 
finitely-generated abelian group $\Z^n$. The matrix $M$ defines an
endomorphism of the $n$-dimensional torus. The mapping torus of this
endomorphism has fundamental group $\Gamma_M$, and is the phase space of
the suspension semiflow of the endomorphism, a semiflow with partially
hyperbolic dynamics (when $M$ is an automorphism, and so $\Gamma_M$ is
polycyclic, the suspension semiflow is actually a flow). Here we see an
example of how the geometric model of a solvable group is actually the
phase space of a dynamical system.

But Bieri-Strebel \cite{BS1} have shown that {\em every} finitely
presented solvable group is, up to commensurability, an ascending HNN
extension with base group a finitely generated solvable group.  In this
way every finitely presented solvable group is the phase space of a
dynamical system, probably realizable geometrically as in the
abelian-by-cyclic case.

\bigskip

%

\end{document}